\documentclass[10pt,a4paper]{amsart} 
\usepackage{graphicx,latexsym,amssymb,diagrams}
\theoremstyle{plain} 
\newtheorem{thm}{Theorem}[section] 
\newtheorem{lem}[thm]{Lemma} 
\newtheorem{prop}[thm]{Proposition} 
 
\theoremstyle{defi} 
\newtheorem{defi}[thm]{Definition} 
\theoremstyle{remark} 
\newtheorem{rem}[thm]{Remark}

\newtheorem{example}[thm]{Example} 
 
\numberwithin{equation}{section}
\numberwithin{figure}{section}
\newcommand{\bd}{\begin{description}}   
\newcommand{\ed}{\end{description}} 
\newcommand{\ba}{\begin{array}}      \newcommand{\ea}{\end{array}} 
\newcommand{\bc}{\begin{center}}     \newcommand{\ec}{\end{center}} 
\newcommand{\be}{\begin{enumerate}}  \newcommand{\ee}{\end{enumerate}} 
\newcommand{\beq}{\begin{eqnarray}}  \newcommand{\eeq}{\end{eqnarray}} 
\newcommand{\beQ}{\begin{eqnarray*}} \newcommand{\eeQ}{\end{eqnarray*}} 
\newcommand{\bi}{\begin{itemize}}    \newcommand{\ei}{\end{itemize}}

\newcommand{\ov}{\overline} 
\newcommand{\we}{\wedge} 
\newcommand{\ve}{\varepsilon} 
\newcommand{\Y}{\mathsf{Y}} 
\newcommand{\I}{\mathsf{I}} 
\newcommand{\nbpt}{18}
\newcommand{\figtotext}[3]{\begin{array}{c}\includegraphics{#3}\end{array}}
\newcommand{\double}{\figtotext{\nbpt}{\nbpt}{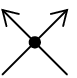}}
\newcommand{\Over}{\figtotext{\nbpt}{\nbpt}{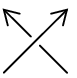}}
\newcommand{\under}{\figtotext{\nbpt}{\nbpt}{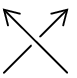}}
\newcommand{\UU}{\figtotext{\nbpt}{\nbpt}{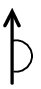}}
\newcommand{\DU}{\figtotext{\nbpt}{\nbpt}{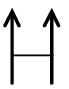}}
\newcommand{\UD}{\figtotext{\nbpt}{\nbpt}{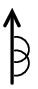}}
\newcommand{\UDi}{\figtotext{\nbpt}{\nbpt}{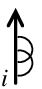}}
\newcommand{\UDB}{\figtotext{\nbpt}{\nbpt}{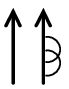}}
\newcommand{\UDC}{\figtotext{\nbpt}{\nbpt}{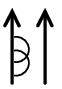}}
\newcommand{\DDA}{\figtotext{\nbpt}{\nbpt}{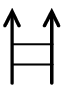}}
\newcommand{\DDAi}{\figtotext{\nbpt}{\nbpt}{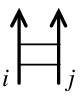}}
\newcommand{\DDBa}{\figtotext{\nbpt}{\nbpt}{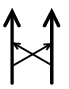}}
\newcommand{\DDCa}{\figtotext{\nbpt}{\nbpt}{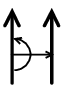}}
\newcommand{\DDDa}{\figtotext{\nbpt}{\nbpt}{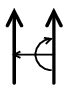}}
\newcommand{\DDB}{\figtotext{\nbpt}{\nbpt}{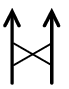}}
\newcommand{\DDBi}{\figtotext{\nbpt}{\nbpt}{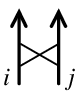}}
\newcommand{\DDC}{\figtotext{\nbpt}{\nbpt}{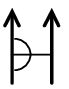}}
\newcommand{\DDCi}{\figtotext{\nbpt}{\nbpt}{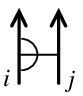}}
\newcommand{\DDCud}{\figtotext{\nbpt}{\nbpt}{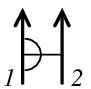}}
\newcommand{\DDCut}{\figtotext{\nbpt}{\nbpt}{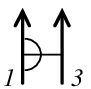}}
\newcommand{\DDCdt}{\figtotext{\nbpt}{\nbpt}{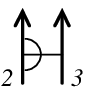}}
\newcommand{\DDD}{\figtotext{\nbpt}{\nbpt}{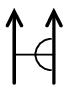}}
\newcommand{\Ddb}{\figtotext{\nbpt}{\nbpt}{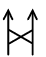}}
\newcommand{\Ddc}{\figtotext{\nbpt}{\nbpt}{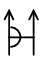}}
\newcommand{\Ddd}{\figtotext{\nbpt}{\nbpt}{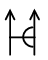}}
\newcommand{\TDA}{\figtotext{\nbpt}{\nbpt}{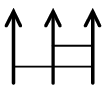}}
\newcommand{\Tda}{\figtotext{\nbpt}{\nbpt}{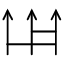}}
\newcommand{\TDB}{\figtotext{\nbpt}{\nbpt}{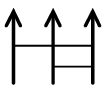}}
\newcommand{\TDC}{\figtotext{\nbpt}{\nbpt}{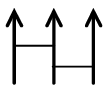}}
\newcommand{\Tdd}{\figtotext{\nbpt}{\nbpt}{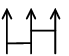}}
\newcommand{\TDD}{\figtotext{\nbpt}{\nbpt}{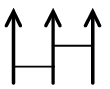}}
\newcommand{\TDE}{\figtotext{\nbpt}{\nbpt}{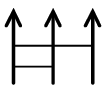}}
\newcommand{\TDF}{\figtotext{\nbpt}{\nbpt}{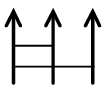}}
\newcommand{\Tde}{\figtotext{\nbpt}{\nbpt}{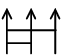}}
\newcommand{\TDAi}{\figtotext{\nbpt}{\nbpt}{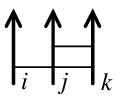}}
\newcommand{\TDBi}{\figtotext{\nbpt}{\nbpt}{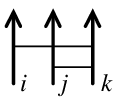}}
\newcommand{\TDCi}{\figtotext{\nbpt}{\nbpt}{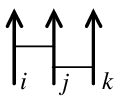}}
\newcommand{\TDDi}{\figtotext{\nbpt}{\nbpt}{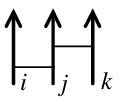}}
\newcommand{\TDEi}{\figtotext{\nbpt}{\nbpt}{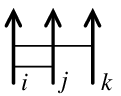}}
\newcommand{\TDFi}{\figtotext{\nbpt}{\nbpt}{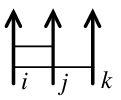}}
\newcommand{\SING}{\figtotext{\nbpt}{\nbpt}{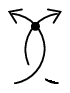}}

\newcommand{\QDAi}{\figtotext{\nbpt}{\nbpt}{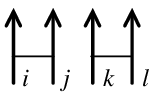}}
\newcommand{\QDBi}{\figtotext{\nbpt}{\nbpt}{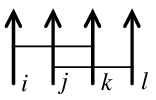}}
\newcommand{\QDDi}{\figtotext{\nbpt}{\nbpt}{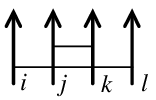}}
\newcommand{\SINGi}{\figtotext{\nbpt}{\nbpt}{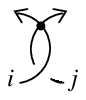}}

\newcommand{\sijk}{\figtotext{\nbpt}{\nbpt}{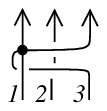}}
\newcommand{\du}{\figtotext{\nbpt}{\nbpt}{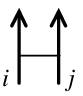}}
\newcommand{\Iij}{\mathsf{I}_{i,j}}
\newcommand{\Iii}{\mathsf{I}_{i,i}}
\newcommand{\Yi}{\mathsf{Y}[(i,1);(i,2);(i,3)]}
\newcommand{\Yii}{\mathsf{Y}[(i,m);(i,n);-]}
\newcommand{\Yiib}{\mathsf{Y}[(i,n);(i,m);-]}
\newcommand{\Yij}{\mathsf{Y}[(i,m);(i,n);j]}
\newcommand{\Yji}{\mathsf{Y}[(j,m);(j,n);i]}
\newcommand{\Yijk}{\mathsf{Y}[i;j;k]}
\newcommand{\Yjik}{\mathsf{Y}[j;i;k]}
\newcommand{\YIJK}{\mathsf{Y}[(i,m);(j,n);(k,l)]}
\newarrow{Noto} ----- 
\newarrow{To} ----> 
\newarrow{Mapsto} |---> 
\newarrow{Into} >---> 
\newarrow{Onto} ----{>>} 
\newarrow{Iso} >---{>>} 
\newarrow{Inclus} C---> 
\newarrow{Equal} ===== 
\newarrow{Dash} {dash}{dash}{dash}{dash}{dash} 
\newarrow{Dashto} {dash}{dash}{dash}{dash}{>} 
\begin{document} 
\title[Vassiliev invariants of order two for string links]{On Vassiliev invariants of order two for string links} 
\author[J.B. Meilhan]{Jean-Baptiste Meilhan} 
\address{Laboratoire de Math\'ematiques Jean Leray\\
         UMR 6629 CNRS/Universit\'e de Nantes \\
         2, rue de la Houssini\`ere \\
         BP 92208 \\
         44322 Nantes Cedex 03\\ 
         France}
	 \email{meilhan@math.univ-nantes.fr}
%\thanks{Commutative diagrams were drawn with Paul Taylor's package}
%\subjclass[2000]{57N10, 57M27}
%\keywords{}
% 
\begin{abstract} 
We show that the Casson knot invariant, linking number and Milnor's triple linking number, together with a 
certain $2$-string link invariant $V_2$, are necessary and sufficient to express any string link Vassiliev 
invariant of order two. Explicit combinatorial formulas are given for these invariants. 
%This result is applied to the study of a conjecture of K. Habiro on $C_k$-equivalence for string links. 
This result is applied to the theory of claspers for string links. 
\end{abstract} 
\maketitle 
\section{Introduction}

\subsection{History and motivations}
Knot theory experienced a major developement in the early 90's through the work of V. Vassiliev \cite{V}, 
restated in simple and combinatorial terms by J. Birman and X.S. Lin \cite{birman,BL}. 
At low degree, Vassiliev invariants are well understood for links : H. Murakami proved that any Vassiliev invariant $v$ of order 
two is a linear combination of linking numbers and their products and Casson invariants of the components, the coefficients 
being given by the initial data of $v$ \cite{mur}. 
When considering the case of string links, which are links with boundary \cite{HL}, more invariants exist at order two : 
Milnor's triple linking numbers are indeed known to be of finite type for string links \cite{BNat,XSL}.
It is therefore a natural question to ask how many (and which) additional invariants are needed to state a Murakami-type result in 
the string link case. \\
Such a result will allow us to give various combinatorial formulas for Vassiliev invariants of order two, in particular for 
Milnor's triple linking numbers. It is also applied to the study of a conjecture of K. Habiro relating Vassiliev invariants to the theory of claspers \cite{habi}.

\subsection{Preliminaries}
Let $D^2$ be the standard two-dimensional disk, and $x_1,...,x_n$ be $n$ marked points in the interior of $D^2$. 
Let us recall from \cite{HL} the definition of a $n$-string link. 
\begin{defi} \label{defsl}
An \emph{$n$-component string link} (or \emph{$n$-string link}) is a proper, smooth embedding 
\[ \sigma : \bigsqcup_{i=1}^n I_i \rTo D^2\times I \]
of $n$ disjoint copies $I_i$ of the unit interval such that, for each $i$, the image $\sigma_i$ of $I_i$ runs from
$(x_i,0)$ to $(x_i,1)$. 
$\sigma_i$ is called the \emph{$i^{th}$ string of $\sigma$}.
\end{defi}
\noindent Note that each string of an $n$-string link is equipped with an (upward) orientation induced by the natural orientation of $I$.\\
The set $\mathcal{SL}(n)$ of all $n$-string links (up to isotopy with respect to the boundary) has a monoid structure, with composition 
given by the \emph{stacking product} and with the trivial $n$-string link $1_n$ as unit element. 

For the study of Vassiliev invariants for string links, we have to consider \emph{singular $n$-string links}, which are 
immersions $\bigsqcup_i I_i \rTo D^2\times I$ whose singularities are transversal double points (in finite number). 
A singular $n$-string link $\sigma$ with $k$ double points can be thought of as an element of 
$\mathbf{Z}\mathcal{SL}(n)$, the free $\mathbf{Z}$-module on $\mathcal{SL}(n)$, by the following well-known skein formula 
\begin{equation} \label{double}
\double = \Over - \under. 
\end{equation}
\begin{defi} \label{defivass}
Let $A$ be an Abelian group. An $n$-string link invariant $f : \mathcal{SL}(n)\rightarrow A$ is a 
\emph{Vassiliev invariant of order $k$} if its natural extension to $\mathbf{Z}\mathcal{SL}(n)$ 
vanishes on every $n$-string-link with (at least) $k+1$ double points.
\end{defi}

A \emph{chord diagram of order $k$} is a disjoint union $\bigsqcup^n_{i=1} I_i$ of $n$ oriented and ordered copies of the unit interval,  
with $k$ chords on it. We denote by $\mathcal{D}^n_k$ the set of chord diagrams of order $k$ (up to order and orientation preserving 
diffeomorphisms of the copies of $I$).
In this paper, the copies of $I$ (called the \emph{strands} of the diagram) will be represented by bold lines, and the chords 
will be drawn with thin lines (dashed lines are also often used in the literature). 

Let $\sigma$ be a singular $n$-string link with $k$ double points. 
The preimages of the double points of $\sigma$ form a subset of $\bigsqcup^n_{i=1} I_i$ which consists in $k$ pairs of points. 
The chord diagram $D$ of order $k$ \emph{associated to $\sigma$} is the element of $\mathcal{D}^n_k$ obtained by joining each of these $k$ pairs 
by a chord ($\sigma$ is said to \emph{respect} the chord diagram $D$). \\
For example, the singular $2$-string link $\SING$ respects the chord diagram $\DU$.

\begin{defi} \label{defpoi}
Let $A$ be an Abelian group. An $A$-valued \emph{weight system of order $k$} is a map 
$W^n_k : \mathcal{D}^n_k\rightarrow A$ satisfying the $(1T)$ and $(4T)$ relations : 
\beQ
(1T) : & W^n_k(\UU)=0 & \\
(4T) : & W^n_k(N)-W^n_k(S)=W^n_k(W)-W^n_k(E)& \\
\eeQ
where $\UU$ is any chord diagram with an isolated chord (\emph{i.e.} disjoint from the other chords) and 
diagrams $N$, $S$, $W$ and $E$ are those of Figure \ref{4T}.\footnote{They are identical outside of a small ball, inside of which 
they are as depicted.}
\begin{figure}[!h]
\begin{center}
\includegraphics{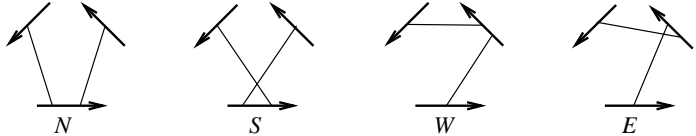}
\caption{The chord diagrams $N$, $S$, $W$ and $E$.} \label{4T}
\end{center}
\end{figure}
\end{defi} 
Note that, for $n=2$ and $3$, the (4T) relation gives : \[ \DDC=\DDA-\DDB=\DDD,\]
\[ \TDA-\TDB=\TDC-\TDD=\TDE-\TDF. \]

Let $v$ be a Vassiliev invariant of order $k$ for $n$-string links. 
For a chord diagram $D$ of order $k$, $v(D)$ denotes the value of $v$ on \emph{any} singular string link respecting $D$. 
It is known (see \cite{birman}) that a Vassiliev invariant satisfies the $(1T)$ and $(4T)$ relations of Definition \ref{defpoi}.
\begin{defi}
The \emph{initial data of $v$} consists in the following :
\begin{displaymath}
\left\{ \begin{array}{ll}
v(D) & \textrm{for every $D\in \mathcal{D}^n_k$.} \\
v(\sigma) & \textrm{for every chord diagram in $\mathcal{D}^n_l$, $l<k$,} \\
 & \textrm{where $\sigma$ is a fixed singular $n$-string link respecting it.} 
\end{array} \right.
\end{displaymath}
\end{defi}
\noindent As outlined in \cite{mur}, two vassiliev invariants of order $k$ with the same initial data coincide. 

For more information about Vassiliev invariants, the reader is refered to \cite{BN}.
\section{Vassiliev invariants of order two for $n$-string links} \label{secVas} 

\subsection{Invariants of order two for $2$-string link} \label{secV2}
\begin{figure}[!h]
\begin{center}
\includegraphics{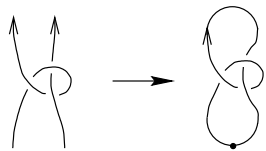}
\caption{The plat-closure of a $2$-string link.} \label{bouclage}
\end{center}
\end{figure}

Recall that the Casson invariant $\varphi(K)$ of a knot $K$ is the coefficient of $z^2$ in the Conway polynomial of $K$ : 
it vanishes on the unknot and equals 1 for the trefoil. Likewise, if $\sigma$ is a $1$-string link, we define its Casson 
invariant by $\varphi(\sigma):=\varphi(\hat{\sigma})$.
Let us denote by $\varphi_i$ the Casson knot invariant of the $i^{th}$ component of a link, 
and by $\mu_{ij}$ the linking number of its $i^{th}$ and $j^{th}$ components.\\
In the case of $2$-component links, it is known \cite{mur} that any Vassiliev invariant $v$ of order two can be written as a linear 
combination 
\begin{equation} \label{sum}
 v=\alpha+\beta.\mu_{12}+\gamma.\varphi_1+\delta.\varphi_2+\omega.\mu^2_{12}, 
\end{equation}
\noindent where $\alpha$, $\beta$, $\gamma$, $\delta$ and $\omega$ are constants. 
When considering the case of string links, it turns out that more invariants are needed :  

Let $V_2:\mathcal{SL}(2)\rTo \mathbf{Z}$ be the map defined by
\begin{equation} \label{eqV2}
V_2(\sigma) : = \varphi\left(p(\sigma)\right) - \varphi(\sigma_1) - \varphi(\sigma_2), 
\end{equation}
where $p(\sigma)$ denotes the plat-closure of $\sigma$ : it is the knot obtained by identifying the origin (resp. the end) of 
$\sigma_1$ with the origin (resp. the end) of $\sigma_2$, with orientation induced by the orientation of $\sigma_1$. 
An example is given in Fig. \ref{bouclage}.\\

As an immediate consequence of the definition of $V_2$, we have :
\begin{prop}
$V_2$ is a $\mathbf{Z}$-valued Vassiliev invariant of order two for $2$-string links.
\end{prop}
\noindent Let us compute this invariant for the string link version $W$ of the Whitehead link, depicted below. 
\bc
\includegraphics{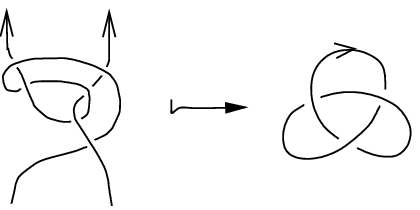}
\ec
The plat-closure of $W$ being the trefoil, and each component being trivial, we have 
\[ V_2(W)=1. \] 
\noindent It follows that $V_2$ can not be expressed as a linear combination of type (\ref{sum}), 
as both the linking number and the Casson invariants of the components vanish on $W$. 
\begin{rem}
Computations (see \cite[\S 2.3.3]{bar}) actually show that there are 4 linearly independent $\mathbf{Z}$-valued Vassiliev invariants 
of order $2$ for $2$-string links. 
This number is equal to the  number of admissible (i.e. without isolated chord) chord diagrams of order $2$ \emph{modulo} 
the $(4T)$ relation : 
\[ \UDC \quad ; \quad \UDB \quad ; \quad \DDA \quad ; \quad \DDC. \]
\noindent The first two correspond to $\varphi_1$ and $\varphi_2$, and the third corresponds to $\mu^2_{12}$. 
The last diagram does not exist in the link case (as it vanishes by $(1T)$ when the endpoints of the strands are identified) : 
we will see in the next section that it actually corresponds to the invariant $V_2$. 
\end{rem}
\begin{rem}
Let $\mathcal{SL}^{0}(2)$ be the submonoid of all 2-string links with vanishing linking number. 
We can prove \cite{these} that the \emph{modulo 2} reduction of $V_2$ coincides on $\mathcal{SL}^{0}(2)$ with the \emph{modulo 2} reduction of 
the Sato-Levine invariant $\beta$ (defined by taking the closure $\hat{\sigma}$ of a 2-string link $\sigma\in \mathcal{SL}^{0}(2)$). 
However, $V_2$ differs from the Sato-Levine invariant. For example, the 2-string link $\sigma$ depicted below is mapped to 0 by $V_2$, 
whereas $\beta(\sigma)=2$.
\bc
\includegraphics{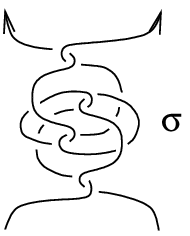}
\ec 
\end{rem}

\subsection{Main result} 
In \cite[Thm. 3.3]{mur}, it is shown that any Vassiliev invariant $v_2$ of order
two for links is of the form 
\[v_2 =a+\sum_{i=1}^{n} b_i.\varphi_i+\sum_{i<j}\left(c_{ij}.\mu_{ij}+d_{ij}.\mu^2_{ij}\right)+
\sum_{\substack{i \\ j<k}}e_{ijk}.\mu_{ij}\mu_{ik}+\sum_{\substack{i<j \\ k<l \\ i<k}}f_{ijkl}.\mu_{ij}\mu_{kl}, \]
where the constants $a$, $b_i$, $c_{ij}$, $d_{ij}$, $e_{ijk}$ and $f_{ijkl}$ are explicitly given in terms of the 
initial data of $v_2$. 
Here, we state a similar result for the string link case. 

Let us denote by $\SINGi$ the singular $n$-string link representing $\du$ depicted below.
\bc
\includegraphics{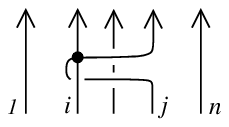}
\ec 
\begin{thm} \label{thmv2}
Let $v$ be a Vassiliev invariant of order two for $n$-string links. 
For any $\sigma\in \mathcal{SL}(n)$, $v(\sigma)$ is given by the following formula \\
\beQ
v(\sigma) & = & v(1_n) + \sum_{i=1}^{n} v(\UDi)\varphi_i(\sigma) \\
	  &+\sum_{i<j}& \left( v(\SINGi)-\frac{1}{2} v(\DDAi) \right) \mu_{ij}(\sigma)+ \\
	  & & \frac{1}{2} v(\DDAi).\mu^2_{ij}(\sigma) - v(\DDCi).(V_2)_{ij}(\sigma) \\
	  &+\sum_{i<j<k}& v(\TDFi).\mu_{ij}\mu_{ik}(\sigma)+v(\TDBi).\mu_{ik}\mu_{jk}(\sigma)+  \\
	  & & v(\TDCi).\mu_{ij}\mu_{jk}(\sigma) + \left( v(\TDCi)-v(\TDDi)\right).\mu_{ijk}(\sigma) \\
   	  &+\sum_{\substack{i<j \\ k<l\textrm{ ; }i<k}}& v(\QDAi).\mu_{ij}\mu_{kl}(\sigma), 
\eeQ
\noindent where $\mu_{ijk}(\sigma)$ denotes Milnor's triple linking number, and $(V_2)_{ij}(\sigma)$ denotes the invariant $V_2$ of 
the $2$-string link $\sigma_i\cup \sigma_j$.
\end{thm}
\noindent \textbf{Proof.} 
Clearly, for arbitrary constants $A$, $B_i$, $C_{ij}$, $D_{ij}$, $E_{ij}$, $F_{ijk}$, $G_{ijk}$, $H_{ijk}$, $P_{ijk}$, $Q_{ijkl}$, 
$R_{ijkl}$ and  $S_{ijkl}$, 
\beQ
w & =& A+\sum_{i=1}^{n} B_i.\varphi_i+\sum_{i<j}\left(C_{ij}.\mu_{ij}+D_{ij}.\mu^2_{ij}+E_{ij}.(V_2)_{ij}\right) \\
 & + & \sum_{i<j<k}\left(F_{ijk}.\mu_{ij}\mu_{ik}+G_{ijk}.\mu_{ik}\mu_{jk}+H_{ijk}.\mu_{ij}\mu_{jk}+P_{ijk}.\mu_{ijk} \right) \\
 & + & \sum_{i<j<k<l}\left(Q_{ijkl}.\mu_{ij}\mu_{kl}+R_{ijkl}.\mu_{ik}\mu_{jl}+S_{ijkl}.\mu_{il}\mu_{jk}\right) 
\eeQ
is a Vassiliev invariant of order two. 
Indeed, recall from \cite{BNat,XSL} that the Milnor $\mu$ invariants of lenght $k$ are Vassiliev invariants of order $k-1$.  
Following the proof of \cite[Thm. 3.3]{mur}, we now choose these constants so that $v$ and $w$ coincide. 

The following computations give us the initial data of $w$ : for all $1\le i<j<k<l\le n$, 
%\[ w(1_n)=A, \]
\beQ
%  & w(1_n)=A,&  \\
w(1_n)=A,\textrm{     } & w(\SINGi)=C_{ij}+D_{ij}, & \\
%w(\UDi)=B_i, & w(\SINGi)=C_{ij}+D_{ij}-P_{ijk}, \\
w(\UDi)=B_i, & w(\DDAi)=2D_{ij},  & w(\DDCi)=-E_{ij}, \\
w(\TDFi)=F_{ijk}, & w(\TDBi)=G_{ijk}, & w(\TDCi)=H_{ijk}, \\
  & w(\TDDi)=H_{ijk}-P_{ijk},& \\
w(\QDAi)=Q_{ijkl}, & w(\QDBi)=R_{ijkl}, & w(\QDDi)=S_{ijkl}. 
\eeQ
Though these computations will not be detailed here, let us develop one of them as an illustration of the general process. 
For example, we consider the chord diagram $\DDCi$ : a singular string ling $\sigma_s$ respecting it is given below. 
\bc
\includegraphics{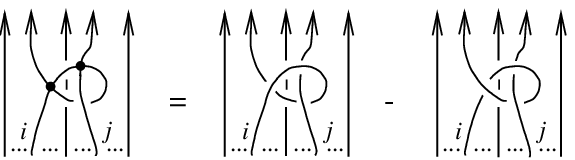}
\ec
By (\ref{double}), $\sigma_s=\sigma_1 - \sigma_2$ where $\sigma_1$ and $\sigma_2$ are the $n$-string links shown in the figure. 
We have $w(\sigma_1)=A+C_{ij}+D_{ij}$ and $w(\sigma_2)=A+C_{ij}+D_{ij}+E_{ij}$, as the plat-closure of $\sigma_1$ 
(resp. $\sigma_2$) is the unknot (resp. the trefoil - see Fig. \ref{bouclage}). The result follows.
$\square$\\
\subsection{Milnor's triple linking number}
In the case of links, the triple linking number $\mu_{ijk}$ is only defined modulo the linking numbers of the components. 
We set $\ov{\mu}_{ijk}\equiv \mu_{ijk}\quad mod \textrm{ gcd($\mu_{ij}$,$\mu_{ik}$,$\mu_{jk}$)}$. $\ov{\mu}_{ijk}$ is known to be 
invariant under cyclic permutations of the indices, and to change its sign under exchange of two indices.\\
As a consequence of Theorem \ref{thmv2}, we have the following formulas for such operations on $\mu_{ijk}$, 
with $i<j<k$ :\\
\bc
\begin{tabular}{|c|ccc|}
\hline 
\textbf{Indices Exchange}  & $\mu_{ijk}$ & = & $-\mu_{jik} + \mu_{ik}\mu_{jk}$ \\
 & & = & $-\mu_{ikj} + \mu_{ij}\mu_{ik} - \mu_{ik}$ \\
 & & = & $-\mu_{kji} + \mu_{ik}\mu_{jk} + \mu_{ij}\mu_{ik} - \mu_{ij}\mu_{jk}$\\[0.1cm]
\hline
\textbf{Cyclic Permutation}  & $\mu_{ijk}$ & = & $\mu_{jki} - \mu_{ij}\mu_{jk} + \mu_{ik}\mu_{jk}$\\
 & & = & $\mu_{kij} - \mu_{ij}\mu_{jk} + \mu_{ij}\mu_{ik} - \mu_{ik}$ \\[0.1cm]
\hline
\end{tabular}
\ec
\noindent Indeed, the singular string link $\sigma_s=\sijk$ satisfies $\mu_{312}(\sigma_s)=1$, $\mu_{132}(\sigma_s)=-1$ and
$\mu_{123}(\sigma_s)=\mu_{213}(\sigma_s)=\mu_{321}(\sigma_s)=\mu_{231}(\sigma_s)=0$. \\

Moreover, we can give a formula for the triple linking numbers in terms of the invariant $V_2$. 
\begin{defi}
Let $\sigma\in \mathcal{SL}(3)$ be a $3$-string link. By stacking above $\sigma$ a negative crossing between the first and third strings 
(such that both strings pass over the second one), and then identifying the endpoints $x_1\times \{0\}$ and $x_1\times \{1\}$, we obtain 
a $2$-string link denoted by $\tilde{\sigma}$ (an example is given below). 
\bc
\includegraphics{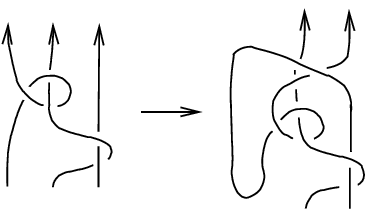}
\ec
\noindent We denote by $V_-$ the Vassiliev invariant of order two $V_-:\mathcal{SL}(3)\rTo \mathbf{Z}$ defined by 
\[ V_-(\sigma):= V_2(\tilde{\sigma}). \]
\end{defi}
%\noindent 
\begin{prop}
The following formula holds for Milnor's triple linking number of $\sigma\in \mathcal{SL}(3)$ :
\begin{equation} \label{eqmu}
\mu_{123}(\sigma) = - V_-(\sigma) + V_2(\sigma_{12}) + V_2(\sigma_{23}) - V_2(\sigma_{13}), 
\end{equation}
where $\sigma_{ij}$ denotes the 2-string link $\sigma_i\cup \sigma_j$ ($1\le i<j\le 3$).
\end{prop}
\noindent \textbf{Proof.} 
Using Theorem \ref{thmv2}, it suffices to compute the initial data of the invariant $V_-$. 
For all $1\le i<j\le 3$, 
\[ V_-(1_3)=0,\quad V_-(\UDi)=0,\quad V_-(\SINGi)=0,\quad V_-(\DDAi)=0. \]
Moreover, 
\[ V_-(\DDCud)=-1,\quad V_-(\DDCut)=-1,\quad V_-(\DDCdt)=+1, \]
and 
\beQ
V_-(\TDB)=0, & V_-(\TDF)=0, & V_-(\TDC)=0, \\
  & V_-(\TDD)=1. &  \square
\eeQ
\begin{rem}
One may ask how revelant is the choice of a negative or positive crossing in the definition of $\tilde{\sigma}$. 
The answer is given in the following formula, which involves the values of $V_-$ and $V_+$ 
(the $3$-string link Vassiliev invariant of order two constructed in a similar way, but with a positive crossing) : 
\begin{equation} \label{pm}
V_-(\sigma) = V_+(\sigma) - \mu_{12}(\sigma). 
\end{equation}
\noindent Indeed, $V_+$ and $V_-$ have the same initial data, except for the singular link $\SINGi$, with $i=1$ and 
$j=2$, on which $V_-$ vanishes, whereas $V_+$ equals $-1$.
\end{rem} 
\section{Combinatorial formulas}
In this section, we give explicit formulas for the invariants of order 2 studied in the previous section, and thus for any 
Vassiliev invariant of order 2 for string links. 
These formulas are derived from a formula of J. Lannes for the Casson knot invariant, given in terms of chord diagrams and weight system. 

\subsection{Lannes formula for the Casson knot invariant.} \label{subc2}
In \cite{La}, a combinatorial formula is given for the Vassiliev knot invariant of order 2, the Casson knot invariant. 
This formula is restated in \cite{tuy} in very explicit terms. Let us briefly recall this formula.
Let $\sigma$ be a $1$-string link (abusing notation, we actually denote by $\sigma$ a fixed diagram of the string link). 
To each crossing $x$ of $\sigma$, we assign two values $\epsilon_x\in \{-1;+1\}$ and $\delta_x\in \{0;1\}$ as follows :
\bi
 \item $\epsilon_x$ is the sign of the crossing : \\
       $\epsilon_x=+1$ if $x=\Over$, and $\epsilon_x=-1$ if $x=\under$,
 \item $\delta_x=0$ if the first branch in $x$ (with respect to the orientation) passes over the second one, and $\delta_x=1$ otherwise.
\ei
Let $W^1_2:\mathcal{D}^1_2\rightarrow \mathbf{Z}$ be the weight system given by $W^1_2(\UD) = 1$. \\
Let $\mathcal{P}_2(\sigma)$ be the set of all (non ordered) pairs of crossings of $\sigma$ : to each element $\{x,y\}$ of $\mathcal{P}_2(\sigma)$, we assign 
an element $D_{x,y}\in \mathcal{D}^1_2$ by replacing $x$ and $y$ by double points and looking at the chord diagram of order two associated to this 
singular string link.\\ 
We have the following formula for the Casson invariant $\varphi$ of a $1$-string link $\sigma$ : 
\[\varphi(\sigma) = \frac{1}{2} \sum_{\{x,y\}\in \mathcal{P}_2(\sigma)} W^1_2(D_{x,y})\epsilon_x \epsilon_y |\delta_x-\delta_y|. \]
\begin{example} Consider a string link version of the trefoil.\\
\parbox[c]{10cm}{
\bc
$\ve_1=\ve_2=\ve_3=-1$, \\ 
$\delta_1=\delta_3=0 \textrm{ and } \delta_2=1$.
\ec}
\parbox[r]{3cm}{\includegraphics{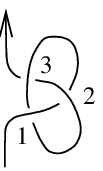}}
So we have $\varphi(\sigma) = 1/2\left( W^1_2(D_{1,2})+W^1_2(D_{2,3})\right)$. 
Now, $D_{1,2}=D_{2,3}=\UD$. \\
It follows that $\varphi(\sigma)=1$.
\end{example}

\subsection{Lannes type formula for Vassiliev invariants of order 2} \label{subv2}
In order to compute the Vassiliev invariants of order 2 listed in Thm. \ref{thmv2}, we will consider string links with $k$ components, 
$1\le k\le 4$. In the case $k=3$, we have several special types of crossings: a crossing between strings 
$\sigma_i$ and $\sigma_j$ ($1\le i\le j\le 3$) is said \emph{of type 1} if $i=1$ and $j=3$, and it is said \emph{of type 2} if 
$i=2$ and $j=3$. 
Moreover, denote by $\sigma_s$ the singular 3-string link $\sijk$. 
Given a 3-string link Vassiliev invariant of order two $v$, we define $\delta_v\in \{0;1\}$ as follows
\begin{displaymath} 
\delta_v = \left\{ \begin{array}{ll}
1 & \textrm{if $v(\sigma_s)$ is not 0,} \\
0 & \textrm{otherwise.}
\end{array} \right.
\end{displaymath}

Now, given an $n$-string link, we can define its curling, which is a $1$-string link :
\begin{defi}
The \emph{curling} of $\sigma \in \mathcal{SL}(n)$, denoted by $C_{\sigma}$, is the element of $\mathcal{SL}(1)$ obtained by attaching, 
for $i\in \{1,...,n-1\}$, the end of the $i^{th}$ string to the origin of the $(i+1)^{th}$, such that the `curls' of $C_{\sigma}$ 
successively pass \emph{over} the first string (producing $n-1$ additionnal positive crossings).
\end{defi}
\begin{figure}[!h]
\begin{center}
\includegraphics{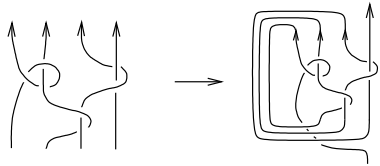}
\caption{The curling of a $4$-string link.} \label{rac}
\end{center}
\end{figure}

Let $\sigma$ be a string link diagram, and $x$ be a crossing of $\sigma$. 
Denote by $\tilde{x}$ the corresponding crossing of the 1-string link $C_\sigma$, and consider the two values 
$\epsilon_{\tilde{x}}\in \{-1;+1\}$ and $\delta_{\tilde{x}}\in \{0;1\}$ associated to $\tilde{x}$, as in \S 3.1.
\begin{prop} \label{lannes}
Any order two Vassiliev invariant $v$ for string links is given by a \emph{Lannes type} formula 
\[ v(\sigma)=\frac{1}{2} \sum_{\{x,y\}\in \mathcal{P}_2(\sigma)} \epsilon_x \epsilon_y |\delta_x-\delta_y|.W_v(D_{xy}), \]
where $W_v$ is a weight system of order two. Here $\epsilon_x$ is the sign of the crossing $x$ and
\begin{displaymath} 
\delta_{x} = \left\{ \begin{array}{ll}
\delta_{v} & \textrm{if $x$ is of type 1,} \\
1-\delta_{\tilde{x}} & \textrm{if $x$ is of type 2,} \\
\delta_{\tilde{x}} & \textrm{otherwise.}
\end{array} \right.
\end{displaymath}
 
\noindent Tables \ref{t1} to \ref{t3} below list the values of the weight system for classical Vassiliev invariants of order 2  
for $n$-string links ($1\le i<j<k<l\le n$).
\end{prop}
\noindent Note that in particular, if $\sigma$ is not a 3-string link, $\delta_{x}=\delta_{\tilde{x}}$ for all crossing $x$ of $\sigma$.

\begin{table}[!h]
\begin{center}
\begin{tabular}{|c|c|ccc|}
\hline
%	            &       &       &        &         \\ 
       v            & $\UDi$ & $\DDAi$ & $\DDBi$ & $\DDCi$  \\
%	            &       &       &        &         \\      
\hline
   $\varphi_{i}$      &   1   &   0   &    0   &   0     \\ 
\hline	
   $\mu_{ij}^2$     &   0   &   2   &    2   &   0     \\ 
   $(V_2)_{ij}$     &   0   &   0   &    1   &  -1     \\
\hline	 
\end{tabular}
\caption{} \label{t1}
\end{center}
\end{table}
\begin{table}[!h]
\begin{center}
\begin{tabular}{|c|cccccc|c|}
\hline
%	            &       &        &        &        &        &        &     	\\ 
       v            & $\TDAi$ & $\TDBi$ & $\TDCi$ & $\TDDi$ & $\TDEi$ & $\TDFi$ \\
%	            &       &        &        &        &        &        &	\\      
\hline 
$\mu_{ij}\mu_{ik}$  &    0    &    0    &    0    &    0    &    1    &   1   \\
$\mu_{ij}\mu_{jk}$  &    0    &    0    &    1    &    1    &    0    &   0   \\
$\mu_{ik}\mu_{jk}$  &    1    &    1    &    0    &    0    &    0    &   0   \\
   $\mu_{ijk}$      &    1    &    0    &    0    &    -1   &    1    &   0   \\
   $\mu_{jik}$      &    0    &    1    &    0    &    1    &    -1   &   0   \\
   $\mu_{kij}$      &    1    &    0    &    1    &    0    &    0    &   -1  \\
   $\mu_{ikj}$      &    -1   &    0    &    0    &    1    &    0    &   1   \\
   $\mu_{jki}$      &    0    &    -1   &    1    &    0    &    1    &   0   \\
   $\mu_{kji}$      &    0    &    1    &    -1   &    0    &    0    &   1   \\
\hline
\end{tabular}
\caption{} \label{t2}
\end{center}
\end{table}
\begin{table}[!h]
\begin{center}
\begin{tabular}{|c|ccc|}
\hline
%	            &       &       &        \\ 
       v            & $\QDAi$  & $\QDBi$ & $\QDDi$  \\
%	            &       &       &            \\      
\hline
$\mu_{ij}\mu_{kl}$  &   1   &   0    &   0    \\
$\mu_{ik}\mu_{jl}$  &   0   &   1    &   0    \\
$\mu_{ik}\mu_{jl}$  &   0   &   0    &   1    \\
\hline	
\end{tabular}
\caption{} \label{t3} 
\end{center}
\end{table}
\noindent \textbf{Proof.} 
The fact that such a Lannes type formula gives rise to a string link invariant essentially follows from the fact that 
a weight system satisfies the $(1T)$ and $(4T)$ relations (Def. \ref{defpoi}). 
Let us now show that such an invariant 
\[ v(\sigma)=\frac{1}{2} \sum_{\{x,y\}\in \mathcal{P}_2(\sigma)} \epsilon_x \epsilon_y |\delta_x-\delta_y|.W_v(D_{xy}), \]
is a Vassiliev invariant of order two. \\
For a pair $\{x,y\}$ of crossings of $\sigma$, we set $\lambda(x,y):=\epsilon_x \epsilon_y |\delta_x-\delta_y|.W_v(D_{xy})$. \\
Let $\sigma_0$ be an $n$-string link with a single double point $x_0$ : by (\ref{double}), 
$\sigma_0=\sigma_+ - \sigma_-$, where $\sigma_+$ (resp. $\sigma_-$) is obtained by replacing $x_0$ by a positive crossing $x_+$ 
(resp. a negative crossing $x_-$). To each crossing $y\ne x_0$ of $\sigma_0$ corresponds a crossing $y\ne x_+$ of $\sigma_+$ and 
a crossing $y\ne x_-$ of $\sigma_-$, such that
\[v(\sigma_{\pm})=\frac{1}{2} \sum_{y\ne x_{\pm}}\lambda(x_{\pm},y)+\frac{1}{2} \sum_{y,z\ne x_{\pm}} \lambda(y,z),\]
\noindent where the notation $y,z\ne x_{\pm}$ stands for a pair of crossings of $\sigma_{\pm}$ 
such that neither $y$ nor $z$ is $x_{\pm}$. 
Now there are two cases. If $x_{\pm}$ is \emph{not} of type 1, we have
\[ \lambda(x_+,y)=\epsilon_{x_+} \epsilon_y |\delta_{x_+}-\delta_y|.W_v(D_{x_+ y}) 
=-\epsilon_{x_-} \epsilon_y |(1-\delta_{x_-})-\delta_y|.W_v(D_{x_- y}).\] 
It follows that $\lambda(x_+,y)-\lambda(x_-,y)=W_v(D_{x_0,y}).\epsilon_y$, for all $y\ne x_0$, and thus 
\beq \label{demi}
v(\sigma_0)=\frac{1}{2} \sum_{y\ne x_0} \epsilon_y.W_v(D_{x_0,y}).
\eeq
If $x_{\pm}$ is of type 1, then 
\[ \lambda(x_+,y)-\lambda(x_-,y)=\epsilon_y |\delta_{v}-\delta_y|.W_v(D_{x_+ y}) 
-(-1) \epsilon_y |\delta_{v}-\delta_y|.W_v(D_{x_- y}). \]
 So in this case 
\beq \label{typ1}
v(\sigma_0)=\sum_{y\ne x_0} \epsilon_y |\delta_{v}-\delta_y|.W_v(D_{x_0,y}).
\eeq
\noindent Similarly, consider a singular string link $\sigma_{12}$ with two double points $x_1$ and $x_2$. 
Suppose first that neither of the two is of type 1. Then, by applying 
(\ref{double}) to, say, $x_2$, we obtain by (\ref{demi}) : 
\beQ
v(\sigma_{12}) & = & \frac{1}{2} W_v(D_{x_1,x_2}) + \frac{1}{2} \sum_{y\ne x_1, x_2} \epsilon_y.W_v(D_{x_1,y}) \\
               & - & \frac{1}{2} W_v(D_{x_1,x_2}).(-1) - \frac{1}{2} \sum_{y\ne x_1, x_2} \epsilon_y.W_v(D_{x_1,y}).
\eeQ
So the value of $v$ on $\sigma_{12}$ only depends on its two double points :
%\beq \label{entier}
$v(\sigma_{12})=W_v(D_{x_1,x_2})$.
%\eeq
%
Now suppose that (say) $x_2$ is of type 1 and $x_1$ is not. 
By applying (\ref{double}) to $x_2$ we likewise obtain $v(\sigma_{12})=W_v(D_{x_1,x_2})$.\\
Finally, suppose that both $x_1$ and $x_2$ are of type 1. Then following (\ref{typ1}) we have: 
\beQ
v(\sigma_{12}) & = & \sum_{y\ne x_1, x_2} \epsilon_y |\delta_{v}-\delta_y|.W_v(D_{x_1,y}) + |\delta_{v}-\delta_v|.W_v(D_{x_1,x_2}) \\
               & - & \sum_{y\ne x_1, x_2} \epsilon_y |\delta_{v}-\delta_y|.W_v(D_{x_1,y}) -  (-1).|\delta_{v}-\delta_v|.W_v(D_{x_1,x_2}),
\eeQ
so in this case, $v(\sigma_{12})=0$.\\
\noindent It follows that, in every cases, $v$ vanishes on any string link with $3$ double points. 

The remaining part of the proof is completed by computing the initial data of such an invariant $v$ for 
each weight system of tables \ref{t1} to \ref{t3}, 
and applying Theorem \ref{thmv2}. These computations are easily performed using the above expressions 
of $v(\sigma_{12})$, for $\sigma_{12}$ a singular string link with two double points.
$\square$\\
\begin{example}
Let $W^2_2:\mathcal{D}^2_2\rightarrow \mathbf{Z}$ be the weight system of order two given by
\bc
\begin{tabular}{c|cccccc}
\hline
        & $\UDB$ & $\UDC$ & $\DDA$ & $\DDB$ & $\DDC$ & $\DDD$  \\
\hline
$W^2_2$ & 0 & 0 & 0 & 1 & -1 & -1 \\
\hline
\end{tabular}
\ec

The $2$-string link invariant $V_2$ is given by the following formula :
\[ V_2(\sigma) = \frac{1}{2} \sum_{\{x,y\}\in \mathcal{P}_2(\sigma)}  W^2_2(D_{x,y})\epsilon_x \epsilon_y |\delta_x-\delta_y|. \]
\noindent In other words, $V_2$ is given by 
\begin{equation} \label{equV2}
V_2(\sigma) = \frac{1}{2} \sum_{\substack{\{x,y\}\in \mathcal{P}_2(\sigma) \\ D_{x,y}=\Ddb}} \epsilon_x \epsilon_y |\delta_x-\delta_y| 
 - \frac{1}{2} \sum_{\substack{\{x,y\}\in \mathcal{P}_2(\sigma) \\ D_{x,y}=\Ddc\textrm{or}\Ddd}} \epsilon_x \epsilon_y |\delta_x-\delta_y|.
\end{equation}

Likewise, we obtain an explicit formula for the triple linking number $\mu_{123}$ :
\begin{equation} \label{equmu3}
\mu_{123}(\sigma) = \frac{1}{2} \sum_{\substack{\{x,y\}\in \mathcal{P}_2(\sigma) \\ D_{x,y}=\Tda\textrm{or}\Tde}}  
\epsilon_x \epsilon_y |\delta_x-\delta_y| 
 - \frac{1}{2} \sum_{\substack{\{x,y\}\in \mathcal{P}_2(\sigma) \\ D_{x,y}=\Tdd}}  
\epsilon_x \epsilon_y |\delta_x-\delta_y|. 
\end{equation} 
\end{example}
\begin{rem}
These Lannes-type formulas are very similar to Gauss diagram formulas developped by T. Fiedler \cite{F} and M. Polyak and O. Viro \cite{pv}. 
In particular, formula (\ref{equmu3}) is to be compared with \cite[Prop. 4.1]{polyak}. Likewise, one has the following Gauss diagram formula 
for the invariant $V_2$, due to M. Polyak \cite{po} :\\

\textit{Let $G_D$ be a Gauss diagram of a 2-string link $\sigma$. Then}
\begin{equation} \label{equGDV2}
V_2(\sigma)=\langle \DDBa - \DDCa - \DDDa, G_D \rangle.
\end{equation}
\end{rem}
\section{$C_k$-equivalence for string links} \label{last}

\subsection{A brief review of clasper theory for string links} \label{review} 
\begin{figure}[!h]
\begin{center}
\includegraphics{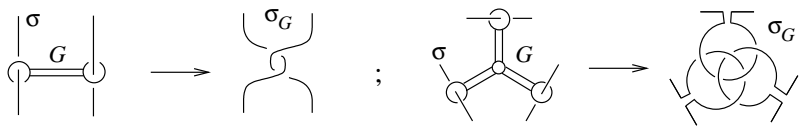}
\caption{A $C_1$-move and a $C_2$-move on a string link $\sigma$.} \label{ckmove}
\end{center}
\end{figure}

Habiro's claspers give a nice reformulation of Vassiliev theory. 
Let us briefly recall from \cite{habi} the basic notions of clasper theory for string links.
\begin{defi}
A \emph{clasper} $G$ for an $n$-string link $\sigma$ in $D^2\times I$ is the embedding 
\[ G: F\rTo D^2\times I \]
of a surface $F$ which is the planar thickening of a uni-trivalent tree.\footnote{A \emph{uni-trivalent tree} is a simply connected 
graph having only univalent and trivalent vertices.} \\
The (thickened) univalent vertices are called the \emph{leaves} of $G$, the trivalent ones 
are the \emph{nodes} of $G$ and we still call \emph{edges} of $G$ the thickened edges of the tree.\\
The string link $\sigma$ is disjoint from $G$, except at the level of each leaf, that it intersects at one point (transversally, in 
the interior).\\
The \emph{degree} of a clasper $G$ is the number of nodes \emph{plus 1}.
\end{defi}
A degree $k$ clasper for a string link $\sigma$ is the instruction for a modification on $\sigma$, called a \emph{$C_k$-move} :
\[ \sigma \mapsto_{C_k} \sigma_G\in \mathcal{SL}(n), \]
which is a surgery move along a framed link defined by $G$. 
More precisely, surgery along the framed link associated to $G$ maps $(D^2\times I)$ to a diffeomorphic $3$-manifold $(D^2\times I)_G$. 
We actually denote by $\sigma_G\subset (D^2\times I)$ the preimage of $\sigma_G\subset (D^2\times I)_G$ by this diffeomorphism. 
Figure \ref{ckmove} shows how a $C_k$-move looks like for $k=1$ or $2$. \\

\noindent Similarly, a $C_k$-move ($k\ge 3$) is realized by a connected sum with a $(k+1)$-component iterated Bing double of the Hopf link. 
Note that a $C_1$-move is equivalent to a crossing change, and that a $C_2$-move is equivalent to a $\Delta$-move \cite{matv,MN}.
\begin{defi}
Let $k\ge 1$ be an integer. The \emph{$C_k$-equivalence}, denoted by $\sim_{C_k}$, is the equivalence relation on string links 
generated by the $C_k$-moves and isotopies with respect to the endpoints.
\end{defi}
\noindent It can be shown that, if $1\le k\le n$, the $C_n$-equivalence relation implies the $C_k$-equivalence.

As said at the begining of the section, clasper theory allows us to reformulate the notion of Vassiliev invariants :\\
\textbf{Definition \ref{defivass}bis.} 
\textit{An $n$-string link invariant $f : \mathcal{SL}(n)\rightarrow A$ is a \emph{Vassiliev in\-va\-riant of order $k$} if, 
for all $n$-string link $\sigma$ and all familly $F=\{G_1,...G_n\}$ of $n$ disjoints claspers for $\sigma$ such that 
$\sum_i deg(G_i)=k+1$, the following equality holds :
\[ \sum_{F'\subseteq F} (-1)^{F'} f\left(\sigma_{F'}\right) = 0. \]}

We conclude this section with some few basic results of calculus of claspers, whose proof can be found in \cite{habi}. 
\begin{prop} \label{facts}
Let $G$ be a degree $k$ clasper for a string link $\sigma$. Then, 
\be
\item[(1)] the $C_{k+1}$-equivalence class of $\sigma_G$ is not modified when an edge of $G$ passes across $\sigma$ or across another 
edge of $G$, or is full-twisted; 
\item[(2)] the $C_{k+1}$-equivalence class of $M_G$ is inverted when an edge of $G$ is half-twisted.
\ee
\end{prop}

\begin{rem}
Every clasper considered here is what is called in \cite{habi} a \emph{simple strict clasper}. 
Note that the definitions and properties given here are also valid for knots and links in $S^3$.
\end{rem}
\subsection{On a conjecture of Habiro on string links} 
String links are extensively studied in \cite[\S 5 and 6]{habi}. The following conjecture
is proposed (see \cite[Conj. 6.13]{habi}).\\
\vspace{0.2cm}
\textbf{Conjecture (Habiro) :} 
\textit{Two $n$-string links in $D^2\times I$ are $C_{k+1}$-equivalent if and only if they are not distinguished by any invariant of 
degree $k$.}
\begin{rem} \label{remconj} 
Note that the \emph{only if} part of the statement is true. 
More precisely, the \emph{only if} part is a general fact that also holds for knot and links in $S^3$. 
The converse statement happens to be true in the knot case \cite[Thm. 6.18]{habi} (see also \cite{funar}), and it is known to fail in the link case 
\cite[Prop. 7.4]{habi}. In this section, we state the \emph{if} part for $k=1$ and $2$. 
\end{rem}
Let us denote by $\mathcal{SL}_{k}(n)\subseteq \mathcal{SL}(n)$ the submonoid of (isotopy classes of) $n$-string links which are 
$C_k$-equivalent to the trivial one. Note that, as a $C_1$-move is just a crossing change, $\mathcal{SL}_{1}(n)=\mathcal{SL}(n)$. 
We have the following descending filtration of monoids 
\[ \mathcal{SL}(n)=\mathcal{SL}_{1}(n)\supset \mathcal{SL}_{2}(n)\supset ... \] 
It is shown in \cite{habi} that 
\[ \overline{\mathcal{SL}}_{k}(n):=\mathcal{SL}_{k}(n)/C_{k+1} \]
is an Abelian group, for all $k\ge 1$. \\
The first of these groups is easily computable : the next theorem claims that it is isomorphic to a space of \emph{struts} (graphs 
made of a single edge) whose vertices are colored by elements of $\{1,...,n\}$. 
More precisely, let us denote by $\Iij$ the strut whose endpoints are colored by two integers $i$ and $j$ : $\mathsf{A}_1(n)$ denotes 
the free Abelian group generated by struts $\Iij$, $i,j\in \{1,...,n\}$, modulo the relation $\Iii=0$.
\begin{thm} \label{thC1} 
There exists a surgery map $\psi_1:\mathsf{A}_1(n)\rTo \ov{\mathcal{SL}}_{1}(n)$ such that $\psi_1$ is an isomorphim of Abelian groups, 
with inverse given by the linking numbers.
\end{thm} 
\noindent The proof of this theorem, together with the definitions of $\psi_1$, is given in \S 4.3. 
This result gives a characterization of $C_2$-equivalence for string links : two string links are $C_2$-equivalent if and only if they 
they have same linking numbers (the Vassiliev invariants of order $1$). 
As a corollary, by considering the closure of string links, it recovers a theorem of H. Murakami and Y. Nakanishi characterizing link-homology 
for links in terms of $\Delta$-moves \cite[Thm. 1.1]{MN}.\\

We now state the main result of this section, which characterize the $C_3$-equivalen\-ce relation for string links in a similar way 
as above : namely, we identify the Abelian group $\ov{\mathcal{SL}}_{2}(n)$ with a certain space of diagrams by means of invariants 
of order two. \\
For that purpose, we define in \S 4.4 the set $\mathsf{A}_2(n)$ generated by $\Y$-shaped diagrams whose univalent vertices are 
colored by elements of $\{1,...,n\}$ and are equipped with a partial order, modulo some relations. 
In \S 4.4, an isomorphism 
\[ \eta : \mathsf{A}_2(n) \rTo \Lambda^{3}H\oplus S^{2}H \]
is defined, where $H$ is the first homology group of 
$D^2\setminus \{x_1,...,x_n\}$, $\Lambda^{3}H$ is the degree three part of the exterior algebra on $H$ and $S^2H$ is the degree two 
part of the symetric algebra. Moreover, 
$\mathsf{A}_2(n)$ produces a combinatorial upper bound for $\ov{\mathcal{SL}}_{2}(n)$, as we define a surjective surgery map 
\[ \mathsf{A}_2(n) \rOnto^{\psi_2} \ov{\mathcal{SL}}_{2}(n). \]
Finally, we will define a homomorphism $(\mu_{3},V_2,\varphi) : \ov{\mathcal{SL}}_{2}(n) \rTo \Lambda^{3}H\oplus S^{2}H$, 
which is a linear combination of the Casson knot invariant, the invariant $V_2$ defined in \S 2.1 and the 
triple linking numbers (and thus is an invariant of order 2), such that the following result holds.
\begin{thm} \label{thC2} 
The following diagram commutes 
\begin{diagram}   
\mathsf{A}_2(n) & \rTo^{\psi_2} &  \ov{\mathcal{SL}}_{2}(n) \\
& \rdTo<{\eta}  & \dTo>{(\mu_{3},V_2,\varphi)} \\
& & \Lambda^{3}H\oplus S^{2}H
\end{diagram} 
with all maps being isomorphisms. 
\end{thm} 
We easily deduce from this theorem that two elements of $\mathcal{SL}_2(n)$ (two $n$-string links with vanishing linking numbers) are $C_3$-equivalent 
if and only if they are not distinguished by Vassiliev invariants of order $2$ (or equivalently, by Theorem \ref{thmv2}, they are not distinguished by 
the Casson knot invariant, nor the invariant $V_2$ and Milnor's triple linking number). This is to be compared with [TY, Thm 1.4] (a clasp-pass 
move is indeed equivalent to a $C_3$-move).
\subsection{Characterization of $C_2$-equivalence : proof of Thm. \ref{thC1}} \label{subC1}
\subsubsection{Combinatorial upper bound for $\overline{\mathcal{SL}}_{1}(n)$}
Let us recall that $\mathsf{A}_1(n)$ is the Abelian group generated by struts $\Iij$, $1\le i\ne j\le n$. 
We define here a \emph{surgery map} 
\[ \psi_1 : \mathsf{A}_1(n) \rightarrow \overline{\mathcal{SL}}_{1}(n) \]
as follows. Let $1_n$ be the trivial $n$-string link in $(D^2\times I)$. 
For each generator $\Iij$ of $\mathsf{A}_1(n)$, consider in $(D^2\times \{1\})\subset \partial (D^2\times I)$ 
the two disks $D_i$ and $D_j$, which are neighbourhoods of the standard points $x_i$ and $x_j$. Their boundary is equipped with a 
natural orientation. Now, push these disks down inside $(D^2\times I)$ along the appropriate strands of $1_n$ (which they always 
intersect transversally at their center), and connect them with a band in a way which is compatible with the orientations. \\
The surface we obtain is a strict and simple \emph{basic clasper} for $1_n$, that we denote by $\phi(\Iij)$ : 
it follows from Prop. \ref{facts}(1) (which essentially implies that the edge of the clasper can be arbitrarily embeded) that the 
$C_2$-equivalence class of $(1_n)_{\phi(\Iij)}$ does not depend on the choice of $\phi$. 
We set
\[ \psi_1(\I_{ij}):=(1_n)_{\phi(\Iij)} \] 
to be (the $C_2$-equivalence class of) the $n$-string link obtained for $1_n$ by surgery along $\phi(\Iij)$. 
The following result follows from the above observation :
\begin{prop}
$\psi_1:\mathsf{A}_1(n)\rightarrow \overline{\mathcal{SL}}_{1}(n)$ is a well defined, surjective map.
\end{prop}
\noindent Note that the surjectivity is just a consequence of the fact that $\overline{\mathcal{SL}}_{1}(n)$ is generated by elements 
of type $(1_n)_C$, where $C$ is a connected basic clasper. Moreover, $\psi_1$ obviously satisfies the relation $\Iii=0$, thanks to 
Prop. \ref{facts}(1).

\subsubsection{Proof of Theorem \ref{thC1}}
\begin{figure}[!h]
\begin{center}
\includegraphics{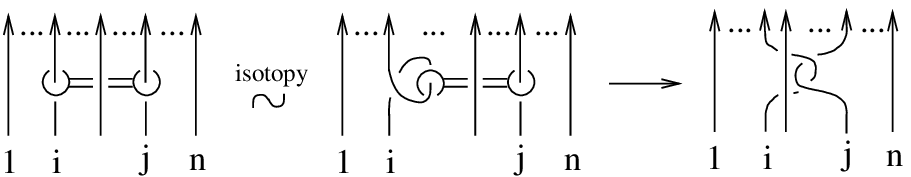}
\caption{} \label{basiq}
\end{center} 
\end{figure}

Let us consider the homomorphism of Abelian groups \\
$\quad \mu_2 : \overline{\mathcal{SL}}_{1}(n) \rightarrow \mathsf{A}_1(n)\quad $  defined by 
\[ \mu_2(\sigma) = - \sum_{1\le i<j\le n} \mu_{ij}(\sigma) \Iij, \]
where $\mu_{ij}(\sigma)$ stands for the linking number of $\sigma_i$ and $\sigma_j$. \\
Given a generator $\Iij$ of $\mathsf{A}_1(n)$, we have just seen that $\psi_1(\Iij)$ is given by the basic clasper depicted in the 
left part of Fig. \ref{basiq}. 
%the $C_2$-equivalence class of the string link $\sigma$ represented below : 

\noindent By a small isotopy on the $i^{th}$ strand, we can use Fig. 4.4 to conclude that $\psi_1(\Iij)$ is 
the $C_2$-equivalence class of the string link depicted in the right part of Fig. \ref{basiq}. 
So we have $\mu_{mn}(\psi_1(\Iij)) = -1$ for $(m,n)=(i,j)$, and $0$ otherwise : $\mu_2\circ \psi_1$ is the identity. 
$\mu_2$ being surjective, it follows that $\psi_1$ is an isomorphism.
\subsection{Characterization of $C_3$-equivalence : proof of Thm. \ref{thC2}} \label{a2}
Let $\mathsf{A}_2(n)$ be the free Abelian group generated by $\Y$-shaped diagrams, whose trivalent vertex is equipped with a cyclic 
ordering on the incident edges and whose univalent vertices are colored by elements of $\{1,...,n\}$, with a total order on vertices 
with the same color, and subject to some relations.  
The notation  
 \[ \Y[n_1;n_2;n_3] \] 
will stand for the $Y$-shaped graph whose univalent vertices are colored by $n_1$, $n_2$ and $n_3\in \{1,...,n\}$ in accordance with the cyclic order, 
(so that our notation is invariant under cyclic permutation of the $n_i$'s). The total order on vertices having the same color is given by a second 
coloring on these vertices, by distinct integers in $\{1;2;3\}$. For example, $\Y[(i,1);(i,2);j]$ and $\Y[(i,2);(i,1);j]$ are two diagrams which only
differ by the order on the $i$-colored vertices ($i\ne j\in \{1,...,n\}$). 
The relations are the following ones:\\[0.1cm]
\beQ
(AS)  : &  \Yijk =  -\Yjik, & \\
(AS2) : &  \Yii = -\Yiib,   & \\
(W)   : &  \Yij =  \Yji,    &
\eeQ \\[0.1cm]
where $i\ne j\ne k\in \{1,...,n\}$ and $m\ne n\in \{1;2;3\}$.
\begin{rem}
Gregor Masbaum pointed out the fact that the Abelian group $\mathsf{A}_2(n)$ actually coincides with the space $\mathcal{B}_2(n)$ 
which consists in $\Y$-shaped and $\Phi$-shaped (2-leg wheels) diagrams modulo the usual AS and IHX relations. Indeed, the $\Y$ part of 
$\mathcal{B}_2(n)$ corresponds to the diagrams of $\mathsf{A}_2(n)$ without repeating colors, and the diagram $\Phi_{ij}$ with univalent vertices colored by $i$ and $j$ 
corresponds to \emph{twice} the diagram $\Yij$ ; $m<n$ (thanks to the Poincar\'e-Birkhoff-Witt isomorphism $\chi$ and the STU relation \cite{BN}). 
Note that the later justifies relation ($W$). 
\end{rem}
\subsubsection{The surgery map $\psi_2$.} 
\begin{figure}[!h]
\begin{center}
\includegraphics{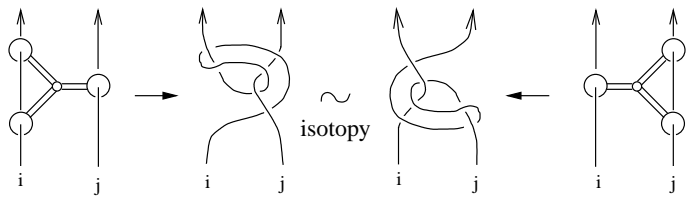}
\caption{The two isotopic string links $w_{ij}$ and $w_{ji}$.} \label{whit}
\end{center}
\end{figure}

Let 
\[ \mathsf{Y}_{i,j,k} := \YIJK \] be a generator of $\mathsf{A}_2(n)$, where $i$, $j$, and $k$ are eventually distinct. 
As in \S 4.3, consider in $(D^2\times \{1\})\subset \partial (D^2\times I)$ the disks $D_i$, $D_j$ and $D_k$, which are 
neighbourhoods of the standard points, and push them inside of $(D^2\times I)$ along the strands of $1_n$, \emph{in the order 
prescribed by the second color}. For example, if $i=j$ and $m>n$, the disk $D_j$ is under $D_i$ in $(D^2\times I)$.
Next, pick an embedded $2$-disk $D$ in the interior of $(D^2\times I)$, disjoint from the $D_i$'s and from $1_n$, orient it in an 
arbitrary way, and connect it to the $D_i$'s with some bands $e_i$ in $(D^2\times I)\setminus 1_n$. These band sums are required to 
be compatible with the orientations, and to be coherent  with the cyclic ordering $(1,2,3)$. The surface we obtain is a degree 2 
clasper for $1_n$ : denote it by $\phi\big(\mathsf{Y}_{i,j,k}\big)$. We denote by 
\[ \psi_2\big(\mathsf{Y}_{i,j,k}\big):=(1_n)_{\phi(\mathsf{Y}_{i,j,k})} \] 
the $n$-string link obtained from $1_n$ by surgery along the clasper $\phi\big(\mathsf{Y}_{i,j,k}\big)$.
\begin{prop} \label{varphi2} 
The $C_3$-equivalence class of $(1_n)_{\phi(\mathsf{Y}_{i,j,k})}$ does not depend on the choice of $\phi$. 
Thus, we have a surjective map : 
\begin{diagram} 
\mathsf{A}_2(n) & \rOnto^{\psi_2} & \overline{\mathcal{SL}}_{2}(n). 
\end{diagram} 
\end{prop} 
\noindent \textbf{Proof.}   
The independence on the choice of the disk $D$, its orientation and the edges $e_i$ is a concequence of Prop. \ref{facts}. \\
We now show that $\psi_2$ satisfies the relations (AS), (AS2) and (W) of $\mathsf{A}_2(n)$. 
The (AS) and (AS2) relations are proved in $\overline{\mathcal{SL}}_{2}(n)$ from Prop. \ref{facts}(2) and an isotopy of the clasper.  
For relation (W), note that a representative for $\psi_2(\Yji)$ is the string link $w_{ij}$ shown in the left part of 
Fig. \ref{whit}, whose $i^{th}$ and $j^{th}$ strings form a Whitehead-type link. As for the link case, this string link version of the 
Whitehead link happens to be symetric (i.e. there is an isotopy which exchanges the role of its components) : $w_{ij}$ is isotopic to the 
string link $w_{ji}$ shown in the right part of the figure, and which is a representative for $\psi_2(\Yij)$.

\noindent The fact that $\psi_2$ is surjective follows from the fact that the Abelian group $\overline{\mathcal{SL}}_{2}(n)$ is generated by 
the $n$-string links $(1_n)_{G}$ where $G$ is a connected degree $2$ clasper. 
$\square$

\subsubsection{The isomorphism $\eta$.} 
Let $\mathsf{A}_{2,k}(n)\subset \mathsf{A}_{2}(n)$ be the subgroup generated by all diagrams whose set of univalent vertices is colored 
with $k$ distinct elements of $\{1,...,n\}$. Note in particular that, as diagrams in $\mathsf{A}_{2,3}(n)$ are colored by 
three distinct elements, there is no second color on their vertices. 
The relations in $\mathsf{A}_{2}(n)$ are graded, so we have 
\begin{equation} \label{truc}
\mathsf{A}_2(n) = \mathsf{A}_{2,3}(n)\oplus \mathsf{A}_{2,2}(n)\oplus \mathsf{A}_{2,1}(n). 
\end{equation} 
Given two distinct integers $m$ and $n$, we denote by $\ve_{mn}\in \{-1;+1\}$ the sign of $(n-m)$. Moreover, 
given $1\le m\ne n\ne l\le 3$, we denote by $\ve_{mnl}$ the signature of the permutation 
$\left[ \begin{array}{ccc}
1 & 2 & 3 \\
m & n & l 
\end{array} \right]$.\\
Consider in $D^2\setminus \{x_1,...,x_n\}$ the $n$ oriented curves defined by a collar of the $n$ boundary components created by  
the removal of the $n$ standard points in $D^2$. We denote by $\{e_1,...,e_n\}$ the associated basis of 
$H=H_1(D^2\setminus \{x_1,...,x_n\};\mathbf{Z})$. \\
Let 
$\quad \eta: \mathsf{A}_2(n) \rightarrow \Lambda^3H\oplus S^2H$ be defined by\footnote{ 
Note that, in this paper, the symetric algebra $SH$ is seen as a \emph{quotient} of the tensor algebra. Therefore, in the following, 
$e_1\otimes e_2\in S^2H$ denotes the class of the tensor product in this quotient.}
\beQ
\eta(\Yijk) = e_i\we e_j\we e_k & \textrm{ on }\mathsf{A}_{2,3}(n), \\
\eta(\Yij) = \ve_{mn}.(e_i\otimes e_j) & \textrm{ on }\mathsf{A}_{2,2}(n), \\
\eta(\Yi) = \ve_{mnl}.(e_i\otimes e_i) & \textrm{ on }\mathsf{A}_{2,1}(n).
\eeQ
\begin{lem}
$\eta$ is an isomorphism of Abelian groups. 
\end{lem}
\noindent \textbf{Proof.} 
It follows from relations (AS), (AS2) and (W) that $\eta$ is well-defined, and thus surjective. Moreover, 
\bi
\item By relation (AS), a basis for $\mathsf{A}_{2,3}(n)$ is given by the set of diagrams of type $\Yijk$ with $i<j<k$. 
	It is therefore a free $\mathbf{Z}$-module of rank $C^n_3$, mapped by $\eta$ onto $\Lambda^3H$ : $\eta$ (retricted to 
	$\mathsf{A}_{2,3}(n)$) defines an epimorphism 
	between two $\mathbf{Z}$-modules of rank $C^n_3$, and then is an isomorphism. 
\item Relations (AS2) and (W) imply that $\mathsf{A}_{2,2}(n)$ (resp. $\mathsf{A}_{2,1}(n)$) has basis given by the 
        $\Yij$'s, with $i<j$ and $n<m$ (resp. the $\Yi$'s). 
	Thus $\mathsf{A}_{2,2}(n)\oplus \mathsf{A}_{2,1}(n)$ is a free $\mathbf{Z}$-module of rank 
	$C^n_2+n=\frac{n(n-1)}{2}+n=\frac{n^2 + n}{2}=rg(S^2H)$, mapped onto $S^2H$ by $\eta$. 
	It follows that $\mathsf{A}_{2,2}(n)\oplus \mathsf{A}_{2,1}(n)\simeq S^2H$ \emph{via} $\eta$. $\square$
\ei
\subsubsection{Proof of Theorem \ref{thC2}.}
Let
$$(\mu_{3},V_2,\varphi) : \mathcal{SL}_{2}(n)\rTo \Lambda^{3}H\oplus S^{2}H$$ 
be the map defined, for all $\sigma\in \mathcal{SL}_{2}(n)$, by
\beQ
(\mu_{3},V_2,\varphi)(\sigma)  & := & \sum_{1\le i<j<k\le n} \mu_{ijk}(\sigma).e_i\we e_j\we e_k \\ 
 & & - \sum_{1\le i<j\le n} V_2(\sigma_i\cup \sigma_j).e_i\otimes e_j  + \sum_{1\le i\le n} \varphi_i(\sigma).e_i\otimes e_i. 
\eeQ
It is a well-defined map which, by Rem. \ref{remconj}, factors to a homomorphism of Abelian groups
\[ (\mu_{3},V_2,\varphi) : \ov{\mathcal{SL}}_{2}(n)\rTo \Lambda^{3}H\oplus S^{2}H. \]
\noindent The following lemma is the last step in proving Theorem \ref{thC2} : it indeed implies that $\psi_2$ and 
$(\mu_{3},V_2,\varphi)$ are isomorphisms. 
\begin{lem} \label{commut}
The following diagram commutes. 
\begin{diagram}   
\mathsf{A}_2(n) & \rTo^{\psi_2} &  \ov{\mathcal{SL}}_{2}(n) \\
& \rdTo<{\eta}  & \dTo>{(\mu_{3},V_2,\varphi)} \\
& & \Lambda^{3}H\oplus S^{2}H
\end{diagram}
\end{lem}
\noindent \textbf{Proof.} 
Let $\Y$ be a generator of $\mathsf{A}_2(n)$ : there are three types of $\Y$ (in the sense of (\ref{truc})). 
We prove that, in these three cases, $(\mu_{3},V_2,\varphi)\circ \psi_2(\Y) = \eta(\Y)$. 

$\bullet$ If $\Y$ is of type $\Yijk$, with $i<j<k$, we have to show that 
$(\mu_{3},V_2,\varphi)\circ \psi_2(\Y) = e_i\we e_j\we e_k \in \Lambda^3H$. \\
A representative for $\psi_2(\Y)\in \ov{\mathcal{SL}}_{2}(n)$ is the string link $\sigma_{ijk}$ obtained from $1_n$ 
by performing a connected sum on strings $\sigma_i$, $\sigma_j$ and $\sigma_k$ with the three components of a Borromean ring (see \S 4.1). 
It follows that $\mu_{abc}(\sigma_{ijk}) = 1$ for $(a,b,c)=(i,j,k)$, and $0$ otherwise. 
Moreover, as every pair of strings of $\sigma_{ijk}$ forms the trivial $2$-string link, $V_2$ is always zero ; the same holds for 
$\varphi$. This proves that $(\mu_{3},V_2,\varphi)\circ \psi_2(\Y) = e_i\we e_j\we e_k$.

$\bullet$ Suppose $\Y$ is of type $\Yij$, with $i<j$ and $m<n$. As seen in the proof of Prop. \ref{varphi2}, a representative 
for $\psi_2(\Y)$ is $w_{ji}$, whose $i^{th}$ and $j^{th}$ strings form a Whitehead-type link (see Fig. \ref{whit}). 
Now, $V_2(\sigma_i\cup \sigma_j)$ equals -1 (as the plat-closure of $w_{ji}$ is the Figure 8 knot), and it vanishes for any other pair of strings. 
All the Milnor's triple linking numbers are zero (as we always are in the case of a $3$-string link with at least one isolated 
string), as well as the Casson knot invariants. It follows that 
$(\mu_{3},V_2,\varphi)\circ \psi_2(\Y) = e_i\otimes e_j = \eta(\Y)\in S^2H$. 

$\bullet$ Finally, let $\Y$ be a generator of type $\Yi$. A representative for its image under $\psi_2$ is $T_{i}$, which 
only differ from $1_n$ by a copy of the trefoil on the $i^{th}$ string (see the figure below).
\bc
\includegraphics{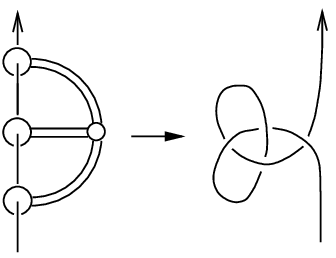}
\ec
\noindent All the triple linking numbers and $V_2$ invariants clearly vanish on $T_i$,  
but the Casson invariant of the $i^{th}$ string equals $1$  
(and $0$ for any other string). \\
So we obtain : $(\mu_{3},V_2,\varphi)\circ \psi_2(\Y) = e_i\otimes e_i =
\eta(\Y)\in S^2H$. \\
This completes the proof of Lemma \ref{commut}, and thus the proof of Theorem \ref{thC2}.
$\square$\\
\\
\\
\textbf{Acknowledgments} 
\\

The author would like to thank Blake Mellor for instructive discussions on Milnor invariants and 
Hitoshi Murakami for stimulating conversations on Vassiliev invariants.  
Acknowlegments are also due to Nathan Habegger and Michael Polyak for valuable comments. 
\\
\\
%\textbf{References} 
\\
 
\vspace{0.5cm}
\noindent
\footnotesize{Commutative diagrams were drawn with Paul Taylor's package.}
\end{document}